\documentclass{article}
\usepackage{amsmath,amssymb,amsthm,amscd}
\title{On the Topology of the Space of\\Contact Structures on Torus Bundles}
\author{Hansj\"org Geiges and Jes\'us Gonzalo}
\date{}

\def\co{\colon\thinspace}
\begin{document}
\newtheorem{thm}{Theorem}
\newtheorem{prop}[thm]{Proposition}
\newtheorem{lem}[thm]{Lemma}
\newtheorem{cor}[thm]{Corollary}
\newtheorem{defn}[thm]{Definition}
\maketitle

\begin{abstract}
We prove the existence of essential loops in the space of contact
structures on torus bundles over the circle.
\end{abstract}

\section{Introduction}
Very little seems to be known about the topology of the space $\Xi (M)$ of
contact structures on a given manifold~$M$. (There is, as yet, no standard
notation for that space. Since $\xi$ is a customary notation for
contact structures, the letter $\Xi$ seems an apt choice.)
Some general results about the
structure of $\Xi$ can be found in \cite[Section~2.4]{hami82} and
\cite{maxi82} (the latter has to be read with a certain amount of caution;
see Mathematical Reviews 83k:58018). For our purposes 
we only need to observe that $\Xi (M)$ is an open subset
(in the $C^k$-topology for $k\geq 1$) of the Grassmannian manifold of
all codimension one subbundles of the tangent bundle~$TM$.

There also appear to be few results concerning $\Xi (M)$ for specific~$M$.
One such result is due to Eliashberg. Let $\xi_0=\ker (x\, dy -y\, dx +z\, dt
-t\, dz)$ be the standard contact structure on $S^3\subset\mathbb{R}^4$
and denote by $\Xi_0(S^3,p_0)$ the subspace of $\Xi (S^3)$ which consists
of contact structures isotopic to $\xi_0$ and with $\xi (p_0)=\xi_0(p_0)$
for some fixed $p_0\in S^3$. In~\cite{elia92} Eliashberg has shown that
$\Xi_0(S^3,p_0)$ is contractible.

This result can be rephrased as follows. Write $\mbox{Diff}^+_{p_0}(S^3)$
for the group of orientation preserving diffeomorphisms of $S^3$ that
fix the contact plane $\xi_0(p_0)$, and denote by $\mbox{Cont}(S^3,\xi_0)$
the subgroup of $\mbox{Diff}^+_{p_0}(S^3)$ consisting of diffeomorphisms
which preserve $\xi_0$.
Then Eliashberg's theorem says that the inclusion of
$\mbox{Cont}(S^3,\xi_0)$ in $\mbox{Diff}^+_{p_0}(S^3)$ is a homotopy
equivalence, cf.~Section~\ref{section:homotopy}. Giroux~\cite{giro01}
gives results about the connected components of other such
contactomorphism groups.

The present paper is concerned with the topology of the space of
contact structures on $T^2$-bundles over~$S^1$.
Specifically, we provide an essentially elementary argument showing
the fundamental group of these spaces (with base point specified below)
to contain an infinite cyclic
subgroup. The basic techniques of the proof are classical
(Gray stability and contact Hamiltonians), but a key ingredient is the 
recent work of Giroux~\cite{giro99} and Kanda~\cite{kand97}
on the classification of contact structures on the $3$-torus.
Throughout we assume contact structures to be coorientable, so they can be 
defined by global contact forms.

For results about related homotopical questions in symplectic
topology see for instance~\cite{abmc99} and~\cite{seid97}.
\section{Contact geometry of the $3$-torus}
Let $(x,y,\theta )$ be standard angular coordinates on $T^3=({\mathbb R}/2\pi
{\mathbb Z})^3$. Define, for $n\in {\mathbb N}$,
\[ \alpha^n_s=\cos (n\theta-2\pi s)\, dx-\sin (n\theta -2\pi s)\, dy,\;\;
s\in [0,1],\]
and $\zeta^n_s=\ker\alpha^n_s$. We shall abbreviate $\zeta^1_s$
to~$\zeta_s$.

Giroux and Kanda have shown that the $\zeta^n_0$ provide a complete list,
up to diffeomorphism, of so-called {\em tight} contact structures on the
3-torus (see~\cite{elia92} for the definition of `tight'). Notice,
however, that the isotopy classification is more subtle. According to
Eliashberg
and Polterovich~\cite{elpo94}, for $\phi\in\mbox{SL}(3,
{\mathbb Z})\subset\mbox{Diff}^+(T^3)$ the contact structures $\phi_*\zeta_0$
and $\zeta_0$ are isotopic if and only if $\phi$ fixes the subspace
${\mathbb Z}^2 \oplus 0$ of the first homology group.
Taken together, these results
yield complete information about the set of connected components of
$\Xi (T^3)$, it being known that $\pi_0(\mbox{Diff}^+(T^3))=\mbox{SL}(3,
{\mathbb Z})$.

The following is our main result about the topology of $\Xi (T^3)$. Below we
shall formulate corresponding statements for all other torus bundles
over~$S^1$.

\begin{prop}
\label{prop:torus}
For each $n\in {\mathbb N}$ the fundamental group $\pi_1(\Xi (T^3),\zeta_0^n)$
based at $\zeta_0^n$ contains an infinite cyclic subgroup, generated by the
loop $\left\{ \zeta_s^n\co 0\leq s\leq 1\right\}$.
\end{prop}

Together with the long homotopy exact sequence of the next section and the
known results about the homotopy type of $\mbox{\rm Diff} (T^3)$ this can be
related to the topology of contactomorphism groups.

Notice that the loop $\{ \xi_s\}$ is defined by a {\em linear} circle
of contact forms,
\[ \xi_s=\ker (\cos (2\pi s)\alpha_0+\sin (2\pi s)\alpha_{1/4}).\]
A pair of contact forms $(\alpha_0,\alpha_{1/4})$
with this property that any non-trivial linear
combination of these forms is again a contact form was called a
{\em contact circle} in \cite{gego95,gego97}, and it is natural to ask
what can be said about such contact circles as elements of $\pi_1(\Xi (M))$.
In \cite{gego97} it is shown that contact circles exist on all closed,
orientable 3-manifolds.

Our proof of Proposition~\ref{prop:torus} makes essential use of the
diffeomorphism classification of contact
structures on $T^3$. It is worth emphasising that the analogue of
Proposition~\ref{prop:torus} for more general torus bundles over the
circle,
formulated in Section~\ref{section:torusbundle} below, rests equally
on the classification of contact structures on~$T^3$, and not
on the classification of contact structures on these spaces themselves.
Again the essential loop is a contact circle.

Clearly, if a contact circle $(\beta_1,\beta_2)$ extends to a {\em contact
sphere} $(\beta_1,\beta_2,\beta_3)$, i.e.\ a triple of contact
forms such that any non-trivial linear combination is a contact form, then
$(\beta_1,\beta_2)$ defines the trivial element in~$\pi_1(\Xi (M))$.

On the other hand, evaluation at a point $p\in M$ and identification of
cooriented radial 2-planes in ${\mathbb R}^3$ with $S^2$ defines a
diffeomorphism
\[ \begin{array}{rcl}
{\mathbb R}^3\supset S^2 & \longrightarrow & S^2\\
(\lambda_1,\lambda_2,\lambda_3) & \longmapsto & \ker (\lambda_1\alpha_1
          +\lambda_2\alpha_2+\lambda_3\alpha_3)(p),
\end{array} \]
so a contact sphere always maps to the generator of $\pi_2(S^2)$ under
this evaluation map and hence defines an element of infinite order
in~$\pi_2(\Xi (M))$. Contact spheres exist, for instance, on the connected sum
of left-quotients of $\mbox{\rm SU}(2)$ and copies of~$S^1\times S^2$,
see~\cite[Prop.~5.7]{gego97}. This does not contradict Eliashberg's result
about the contractibility of $\Xi_0(S^3,p_0)$, of 
course, since in $\Xi (S^3)$ we do not fix the contact plane at a point.
\section{A homotopy exact sequence}
\label{section:homotopy}
The following considerations rest on the well-known concepts of
Gray stability and contact Hamiltonians, cf.~\cite{lima87}, \cite{geig}.
We only recall briefly the facts that we shall use:

\vspace{1mm}

{\bf Gray stability:} Given a smooth family of contact structures $\xi_t$,
$t\in [0,1]$, on a closed manifold~$M$, there is a canonically defined
time-dependent vector field $X_t$ tangent to $\xi_t$ whose flow
$\psi_t$ satisfies $\psi_{t*}\xi_0=\xi_t$ (that is, the differential
$T\psi_t$ maps $\xi_0$ to~$\xi_t$). Given a family of contact forms
$\alpha_t$ defining $\xi_t$, with Reeb vector field~$R_t$,
this vector field is determined by the equations
$\alpha_t(X_t)=0$ and $i_{X_t}d\alpha_t=\mu_t\alpha_t-\dot{\alpha}_t$,
with the function $\mu_t$ determined by $\mu_t=\dot{\alpha}_t(R_t)$, and
with $\dot{\alpha}_t$ denoting the time derivative of~$\alpha_t$.
This $X_t$ is independent of the choice of~$\alpha_t$, even though $\mu_t$
does depend on that choice.

\vspace{1mm}

{\bf Contact Hamiltonians:}
Given a smooth family of functions $H_t$ on a contact manifold $(M,\xi=\ker
\alpha )$, a time-dependent vector field $X_t$ whose flow preserves $\xi$
is defined by
\[ X_t=H_tR+Y_t,\]
where $R$ is the Reeb vector field of $\alpha$ and $Y_t$ is determined by
\[ \alpha (Y_t)=0\]
and
\[ i_{Y_t}d\alpha =dH_t(R)\alpha -dH_t.\]
Indeed, the infinitesimal condition $L_{X_t}\alpha
=\rho_t\alpha$ is equivalent to
the equations above if we require $H_t=\alpha (X_t)$. Put differently,
a time-dependent vector field $X_t$ whose flow preserves $\xi=\ker\alpha$
is completely determined by $H_t=\alpha (X_t)$, the {\em Hamiltonian
function} of~$X_t$.

\vspace{2mm}

Let $(M,\xi_0)$ be a closed contact manifold. Write $\mbox{Diff}_0(M)$ for
the identity component of the diffeomorphism group of $M$ and
$\mbox{Cont}_0(M,\xi_0)$ for its subgroup of contactomorphisms, i.e.
\[ \mbox{Cont}_0(M,\xi_0)=\{ \phi\in\mbox{Diff}_0(M)\co \phi_*\xi_0=\xi_0\} .\]
We have an obvious inclusion $i\co\mbox{Cont}_0(M,\xi_0)\rightarrow
\mbox{Diff}_0(M)$ and, by Gray stability, a surjection
\[ \begin{array}{rcl}
\sigma\co \mbox{Diff}_0(M) & \longrightarrow & \Xi_0(M)\\
\phi & \longmapsto & \phi_*\xi_0,
\end{array} \]
where $\Xi_0(M)$ denotes the component of $\Xi (M)$ containing $\xi_0$.

The following is a well-known (folklore) result:

\begin{prop}
\label{prop:homotopy}
The map $\sigma$ has the homotopy lifting property with respect to cubes
$I^n=[0,1]^n$ and smooth homotopies. Hence
there is an exact sequence of homotopy groups
\[... \stackrel{\Delta}{\longrightarrow} \pi_k(\mbox{\rm Cont}_0)
\stackrel{i_{\#}}{\longrightarrow} \pi_k(\mbox{\rm Diff}_0)
\stackrel{\sigma_{\#}}{\longrightarrow} \pi_k(\Xi_0)
\stackrel{\Delta}{\longrightarrow} \pi_{k-1}(\mbox{\rm Cont}_0)
\stackrel{i_{\#}}{\longrightarrow} ...\]
\end{prop}

\noindent {\em Proof.}
We only need to prove the first statement.
The usual proof that a Serre
fibration yields an exact sequence of homotopy groups then still applies to
give the desired conclusion.

This means we are given a commutative diagram
\[ \begin{CD}
I^n\times 0 @>>> \mbox{\rm Diff}_0(M)\\
@VVV             @VV{\sigma}V\\
I^n\times I @>>> \Xi_0(M),
\end{CD} \]
with the bottom map smooth in the $I$-factor; more precisely,
a family $\xi_{s,t}$ of contact structures on $M$,
continuous in $s\in I^n$, smooth in $t\in I$, and with all $t$-derivatives
continuous in~$s$. We also have a family 
$\phi_{s,0}\in\mbox{Diff}_0(M)$ with $(\phi_{s,0})_*\xi_0=\xi_{s,0}$, and the
aim is to find a lifting $I^n\times I\rightarrow\mbox{\rm Diff}_0(M)$.

By Gray stability (with an added parameter~$s$) one finds a family
of vector fields $X_{s,t}$,
continuous in $s$ and smooth in $t$, such that the flow of the $t$-dependent
vector field $X_{s,t}$ defines a family of diffeotopies $\psi_{s,t}$ which
is continuous in $s$ and satisfies $(\psi_{s,t})_*\xi_{s,0}=\xi_{s,t}$.
Set $\phi_{s,t}=\psi_{s,t}\circ\phi_{s,0}$. Then $(\phi_{s,t})_*\xi_0
=\xi_{s,t}$, so $\phi_{s,t}$ is the desired lift 
of~$\xi_{s,t}$. \hfill $\Box$
\section{Proof of Proposition~1}

We only show that $\{\zeta_s=\zeta^1_s\}$ is non-zero in $\pi_1$; a completely
analogous argument applies to any multiple of this loop and to general~$n$.

Our proof is by contradiction. Assume $\{\zeta_s\}$ defines the trivial
element in $\pi_1(\Xi(T^3),\zeta_0)$. Write $\pi$ for the canonical
submersion
\[ T^2\times {\mathbb R}\longrightarrow T^2\times {\mathbb R}/2\pi {\mathbb Z}
\equiv T^3\]
and set $\widetilde{\zeta}_s=\pi^*\zeta_s$ for each $s\in [0,1]$.
Then $\{\widetilde{\zeta}_s\}$ would define
the trivial element in $\pi_1(\Xi_0^{\mathbb Z}(T^2\times {\mathbb R}))$,
where $\Xi_0^{\mathbb Z}(T^2\times {\mathbb R})$ denotes the connected
component containing $\widetilde{\zeta}_0$ of the space of contact structures
on $T^2\times {\mathbb R}$ which are invariant under shifts by $2\pi$ in
${\mathbb R}$-direction.

Proposition~\ref{prop:homotopy}
applies to $T^2\times {\mathbb R}$ if we replace $\Xi_0$ by
$\Xi_0^{\mathbb Z}(T^2\times {\mathbb R})$ and $\mbox{Diff}_0$ by 
$\mbox{Diff}_0^{\mathbb Z}(T^2\times {\mathbb R})$, the identity component
of the ${\mathbb Z}$-equivariant diffeomorphisms
of $T^2\times {\mathbb R}$, i.e.
\[ \mbox{Diff}_0^{\mathbb Z}(T^2\times {\mathbb R})=\left\{ \phi\in
\mbox{Diff}(T^2\times {\mathbb R})\co \phi (x,y,\theta+2\pi )=\phi (x,y,\theta )
+(0,0,2\pi )\right\}_0,\]
for the ${\mathbb Z}$-invariance resp.\ -equivariance guarantees that the
vector field $X_{s,t}$ used in the proof of
Proposition~\ref{prop:homotopy} still
integrates to a global flow.

We would then have a homotopy $\{\widetilde{\zeta}_{s,t}\}$
rel~$\{ 0,1\}$ of $\{ \widetilde{\zeta}_s\}$
to the constant loop $\{\widetilde{\zeta}_0\}$. We may
assume that $\widetilde{\zeta}_{s,t}$ is smooth in $t$ (by smoothing a homotopy
which is constant near $t=0,1$).

Define a diffeomorphism
\[ \begin{array}{rcl}
\widetilde{\phi }_s\co T^2\times {\mathbb R} & \longrightarrow &
T^2\times {\mathbb R}\\
(x,y,\theta ) & \longmapsto & (x,y,\theta +2\pi s ),
\end{array} \]
so that $\widetilde{\phi}_{s*}\widetilde{\zeta}_0=\widetilde{\zeta}_{s,0}
=\widetilde{\zeta}_0$. By Proposition~\ref{prop:homotopy} we find
a path $\psi(s)$ in
$\mbox{Cont}_0(T^2\times {\mathbb R},\widetilde{\zeta}_0)$
that joins $\widetilde{\phi}_0$ with $\widetilde{\phi}_1$.

Let $X_s$ be the time-dependent vector field that integrates to
$\psi (s)$ at time~$s$.
Choose a contact form $\widetilde{\alpha}_0$ representing $\widetilde{\zeta}_0$
and set $H_s=\widetilde{\alpha}_0(X_s)$. Since the flow of
$X_s$ preserves $\widetilde{\zeta}_0$, the Hamiltonian vector field
of $H_s$ coincides with $X_s$.
The time-one-map of this Hamiltonian flow is $\widetilde{\phi}_1$.
During the diffeotopy, the image 
of $T^2\times [0,2\pi ]$ under $\psi (s)$ stays above $T^2\times
\{ -2\pi (k-1)\}$ for some $k\in {\mathbb N}$ sufficiently large.

Now let $\chi (\theta )$ be a smooth function with $\chi (\theta )\equiv 0$
for $\theta\leq -2\pi k$ and $\chi (\theta )\equiv 1$ for
$\theta\geq -2\pi (k-1)$, and set $\overline{H}_s(x,y,\theta )=
H_s(x,y,\theta )\chi (\theta )$. The Hamiltonian vector field
$\overline{X}_s$ of $\overline{H}_s$ still integrates to a global
flow, and the time-one-map of this flow defines a 
$\widetilde{\zeta}_0$-preserving diffeomorphism
\[ T^2\times [-2\pi k,0]\longrightarrow T^2\times [-2\pi k,2\pi ].\]
But this would imply that the contact structures $\zeta^k$
and $\zeta^{k+1}$ on $T^3$
are diffeomorphic, contradicting the work of Giroux and Kanda.
This contradiction proves the proposition. \hfill $\Box$

\vspace{2mm}

Here is an alternative proof of Proposition~\ref{prop:torus},
based on the work of
Eliashberg-Hofer-Salamon~\cite{ehs95}. Again arguing by contradiction, we
assume that $\{ \zeta_s\}$ is trivial in $\pi_1(\Xi(T^3),\zeta_0)$.
Consider the double cover
\[ T^3_{4\pi}\equiv T^2\times {\mathbb R}/4\pi {\mathbb Z}\longrightarrow
T^2\times {\mathbb R}/2\pi {\mathbb Z}\equiv T^3,\]
and write $\zeta_s'$ for the lift of $\zeta_s$. Then $\{\zeta_s'\}$ is also
a contractible loop in the corresponding space of contact structures.

The same argument as before yields a contact isotopy of $(T^3_{4\pi},
\zeta_0')$, beginning at the identity and ending with a shift by $2\pi$ in
$\theta$-direction. In particular, this contact isotopy would
move the pre-Lagrangian torus (see~\cite{ehs95})
$T^2\times\{ 0\}$ to $T^2\times\{ 2\pi\}$ and therefore separate it completely
from itself. This would contradict \cite[Thm.~3.8.3]{ehs95}, which
uses a Floer homology argument to provide a lower bound for the number of
intersection points between a pre-Lagrangian submanifold $\Lambda_0$ and
the image $\Lambda_1'$ of a Legendrian submanifold $\Lambda_1\subset
\Lambda_0$ under a contact isotopy (The conditions $O_1$ and $O_2$
of that theorem are satisfied in the present context).
\section{Contact circles on torus bundles}
\label{section:torusbundle}
We now extend Proposition~\ref{prop:torus}
to all torus bundles over the circle.
For $A\in\mbox{\rm SL}_2{\mathbb Z}$, write $M_A$ for the quotient
of $T^2\times {\mathbb R}$ under the transformation
\[ \widetilde{A}\co
\left( \left( \begin{array}{cc} x\\y\end{array} \right) ,\theta \right)
\longmapsto \left( A\left(\begin{array}{cc} x\\y\end{array} \right) ,
\theta+2\pi \right).\]
On each $M_A$ there is a family of contact structures $\zeta^n$,
$n\in {\mathbb N}$, characterised up to fibre preserving isotopy
by the property that it descends from
a contact structure on ${\mathbb R}^2\times {\mathbb R}$ of the form
\[ \cos f(\theta )\, dx-\sin f(\theta )\, dy=0,\]
invariant under the transformation~$\widetilde{A}$, and with
\[ 2(n-1)\pi < \sup_{\theta\in{\mathbb R}} (f(\theta +1)-f(\theta ))
\leq 2n\pi ,\]
see~\cite{giro99}.

\begin{prop}
For each $A\in\mbox{\rm SL}_2{\mathbb Z}$ and $n\in {\mathbb N}$, the
fundamental group of $\Xi (M_A)$ based at $\zeta^n$ contains an element of
infinite order.
\end{prop}

To prove this proposition it is convenient to use the Thurston geometries
adapted to the different choices of~$A$, as done previously
in~\cite{gego97}. Rather than giving a case by case proof, we illustrate
the method by a particular example which contains all the ideas required
for the general argument.

Assume that $\mbox{\rm trace}\, A\geq 3$. Then $A$ is conjugate
in $\mbox{\rm GL}_2{\mathbb R}$ to a matrix of the form
\[ A'= \left(\begin{array}{cc}
e^{\gamma}& 0\\
0 & e^{-\gamma}
\end{array} \right) \]
and $M_A$ is modelled on the solvable Lie group $\mbox{\rm Sol}^3$ (the
inhomogeneous Lorentz group), cf.~\cite{gego97}. This means that we can
write $M_A$ as a quotient of ${\mathbb R}^3$ under a group of transformations
$\Gamma$ generated by
\begin{eqnarray*}
(x,y,\theta ) & \longmapsto & (x+a_1, y+b_1, \theta )\\
(x,y,\theta ) & \longmapsto & (x+a_2, y+b_2, \theta )\\
(x,y,\theta ) & \longmapsto & (e^{\gamma} x,e^{-\gamma}y, \theta +\gamma ),
\end{eqnarray*}
where the vectors
\[ \left( \begin{array}{cc} a_1\\b_1\end{array} \right) ,\:
\left( \begin{array}{cc} a_2\\b_2\end{array} \right) \]
generate a lattice $\Lambda$ in ${\mathbb R}^2$ invariant under~$A'$.

Now define, for $n\in {\mathbb N}$,
\[ \alpha_s^n=\cos \left(\frac{2\pi n\theta}{\gamma} -2\pi s\right)
e^{-\theta }dx-
\sin\left(\frac{2\pi n\theta}{\gamma}-2\pi s\right)
e^{\theta} dy,\;\; s\in [0,1],\]
and set $\zeta^n_s=\ker\alpha_s^n$. This defines a contact structure on
${\mathbb R}^3$ that descends to~$M_A$. Up to fibre preserving isotopy,
$\zeta_0^n$ is equal to the $\zeta^n$ studied by Giroux. Write
$\widetilde{\zeta}_0$ for the lift of $\zeta_0=\zeta_0^1$ to
$T^2\times {\mathbb R}$.

Denote by $\Xi_0^{A'}(T^2\times {\mathbb R})$ the connected component
of $\widetilde{\zeta}_0$ in the space of contact structures on $T^2\times
{\mathbb R}$ invariant under the transformation
\[ \widetilde{A}'\co (x,y,\theta )\longmapsto (e^{\gamma} x,
e^{-\gamma}y,\theta +\gamma ),\]
where we think of $T^2$ as the quotient of ${\mathbb R}^2$ under the
lattice~$\Lambda$.
Write $\mbox{\rm Diff}_0^{A'}(T^2\times {\mathbb R})$ for the identity
component of the diffeomorphisms of $T^2\times {\mathbb R}$ that
commute with~$\widetilde{A}'$.

Then the argument in the proof of Proposition~\ref{prop:torus}
goes through as before,
with $\Xi_0^{\mathbb Z}(T^2\times{\mathbb R})$ replaced by
$\Xi_0^{A'}(T^2\times{\mathbb R})$, and $\mbox{\rm Diff}_0^{\mathbb Z}
(T^2\times{\mathbb R})$ by $\mbox{\rm Diff}_0^{A'}(T^2\times{\mathbb R})$.
Assuming that $\{\zeta_s\}$ was trivial in $\pi_1(\Xi(M_A),\zeta_0)$, we
obtain a $\widetilde{\zeta}_0$-preserving diffeomorphism
\[ T^2\times [-k\gamma,0]\longrightarrow
T^2\times [-k\gamma ,\gamma ].\]
On $T^2\times \{ l\gamma\}$, $l\in{\mathbb Z}$, the contact form
$\alpha_0^1$ restricts to $e^{-l\gamma}dx$, so the characteristic
foliation on these tori is always given by $dx=0$. Since the characteristic
foliation on a surface determines the germ of the contact structure along that
surface, we can use the identity map on $T^2$ to
glue the ends of $T^2\times [-k\gamma , \varepsilon\gamma ]$,
with $\varepsilon\in\{ 0,1\}$, and obtain a contact structure on $T^3$,
diffeomorphic to the standard structure $\zeta^k$ of $T^3$ for
$\varepsilon =0$, and $\zeta^{k+1}$ for $\varepsilon =1$. The diffeomorphism
above would induce a diffeomorphism between $\zeta^k$ and $\zeta^{k+1}$
on~$T^3$, contradicting once again the work of Giroux and Kanda.

\vspace{2mm}
 
\noindent
\begin{minipage}[t]{6cm}
Hansj\"org Geiges
 
Mathematisches Institut
 
Universit\"at zu K\"oln
 
Weyertal 86--90

50931 K\"oln

Germany
 
e-mail: geiges@math.uni-koeln.de
\end{minipage}
\begin{minipage}[t]{6cm}
Jes\'us Gonzalo
 
Departamento de Matem\'aticas
 
Universidad Aut\'onoma de Madrid
 
28049 Madrid
 
Spain

e-mail: jesus.gonzalo@uam.es
\end{minipage}
\end{document}